\documentclass{amsart}
\usepackage{latexsym}
\usepackage{amsmath}
\usepackage{amssymb}
\DeclareMathOperator{\sgn}{sign}
\DeclareMathOperator{\sat}{sat}
\DeclareMathOperator{\tor}{Tor}

\begin{document}


{\theoremstyle{plain}%
\newtheorem{theorem}{Theorem}[section]  \newtheorem{corollary}[theorem]{Corollary}
  \newtheorem{proposition}[theorem]{Proposition}
  \newtheorem{lemma}[theorem]{Lemma}
  \newtheorem{question}[theorem]{Question}
  \newtheorem{conjecture}[theorem]{Conjecture}
}
{\theoremstyle{remark}
\newtheorem{fact}{Fact}
\newtheorem{remark}[theorem]{Remark}
}

{\theoremstyle{definition}
\newtheorem{definition}[theorem]{Definition}
\newtheorem{example}[theorem]{Example}
}

\newcommand{\bt}{\underline{t}}
\newcommand{\bj}{\underline{j}}
\newcommand{\bi}{\underline{i}}
\newcommand{\dep}{\operatorname{depth}}
\newcommand{\hit}{\operatorname{ht}}
\newcommand{\op}{\overline{P}}
\newcommand{\xab}{X^{\underline{a}}Y^{\underline{b}}}
\newcommand{\Iz}{I_{Z}}
\newcommand{\Ixp}{I_{\xp}}
\newcommand{\Z}{\mathbb{Z}}
\newcommand{\ax}{\alpha_{\X}}
\newcommand{\bx}{\beta_{\X}}
\newcommand{\kxo}{k[x_0,\ldots,x_n]}
\newcommand{\kx}{k[x_1,\ldots,x_n]}
\newcommand{\popo}{\mathbb{P}^1 \times \mathbb{P}^1}
\newcommand{\pr}{\mathbb{P}}
\newcommand{\pn}{\mathbb{P}^n}
\newcommand{\pnpm}{\mathbb{P}^n \times \mathbb{P}^m}
\newcommand{\pnk}{\mathbb{P}^{n_1} \times \cdots \times \mathbb{P}^{n_k}}
\newcommand{\X}{\mathbb{X}}
\newcommand{\A}{\mathbb{A}}
\newcommand{\Y}{\mathbb{Y}}
\newcommand{\N}{\mathbb{N}}
\newcommand{\M}{\mathbb{M}}
\newcommand{\Q}{\mathbb{Q}}
\newcommand{\Ix}{I_{\X}}
\newcommand{\pix}{\pi_1(\X)}
\newcommand{\pixt}{\pi_2(\X)}
\newcommand{\pipi}{\pi_1^{-1}(P_i)}
\newcommand{\piqi}{\pi_2^{-1}(Q_i)}
\newcommand{\qpi}{Q_{P_i}}
\newcommand{\pqi}{P_{Q_i}}
\newcommand{\pitk}{\pi_{2,\ldots,k}}
\newcommand{\prdim}{\operatorname{proj-}\dim}
\newcommand{\kdim}{\operatorname{K-}\dim}
\newcommand{\reg}{\operatorname{reg}}
\newcommand{\regI}{\operatorname{reg({\it I}_\X)}}
\newcommand{\ri}{\operatorname{ri}}
\newcommand{\Hx}{\mathcal{H}_{\X}}
\newcommand{\HH}{\mathcal{H}}


\title{Multigraded regularity: syzygies and fat points}
\thanks{Revised Version: July 21, 2005}
\author{Jessica Sidman}
\address{Department of Mathematics and Statistics \\
451A Clapp Lab\\
Mount Holyoke College\\
South Hadley, MA 01075, USA}
\email{jsidman@mtholyoke.edu}
\author{Adam Van Tuyl}
\address{Department of Mathematical Sciences \\
Lakehead University \\
Thunder Bay, ON P7B 5E1, Canada}
\email{avantuyl@sleet.lakeheadu.ca}


\begin{abstract}  
The Castelnuovo-Mumford regularity
of a graded ring is an important invariant in 
computational commutative algebra, and there is
increasing interest in  multigraded generalizations.  We study connections 
between two recent definitions of multigraded regularity with a view towards
a better understanding of the multigraded Hilbert function of fat point
schemes in $\pnk.$
\end{abstract}

\maketitle


\section*{Introduction}
Let ${\bf k}$ be an algebraically closed field of characteristic zero.  If $M$ is a finitely generated
graded module over a $\Z$-graded polynomial ring over ${\bf k}$, its Castelnuovo-Mumford regularity, denoted $\reg(M),$ is
an invariant that measures the difficulty of computations involving $M.$ 
Recently, several authors (cf. \cite{ACD,HW,MS}) 
have proposed extensions of the
notion of regularity to a multigraded context.

Taking our cue from the study of the Hilbert functions  of fat points
in $\pr^n$ (cf. \cite{CTV,FL,TV}), we apply these new notions of multigraded regularity 
to study the coordinate ring of a scheme of fat points $Z \subseteq 
\pnk$ with the goal of understanding both the nature of regularity in 
a multigraded setting and what regularity may tell us about the coordinate
ring of $Z.$  This paper also complements the investigation 
in \cite{HVT} of the Castelnuovo-Mumford
regularity of fat points in multiprojective spaces.

In the study of the coordinate ring of a fat point scheme in 
$\pr^n$  many
authors have found beautiful relationships between algebra, geometry,
 and combinatorics
(cf. \cite{H} for a survey when $n=2$).  Extensions and generalizations of such
 results to the multigraded setting are potentially of both theoretical
and practical interest.  Schemes of fat points in products of projective spaces
 arise in algebraic geometry in connection with secant varieties of Segre 
varieties (cf. \cite{CGG,CGG2}).  More generally, the base points of
rational maps betweeen higher dimensional varieties may be non-reduced
schemes of points, and in the case of maps between certain surfaces, 
the regularity of the ideals that arise may have implications for the 
implicitization problem in computer-aided design (cf. \cite{CoxII,ZSCC}).

   It is well known that the
 Castelnuovo-Mumford regularity of a finitely generated $\Z$-graded module 
$M$ can be defined 
either in terms of degree bounds for the generators of the syzygy modules
of $M$ or in terms of the vanishing of graded pieces of local 
cohomology modules.  (See \cite{EG}.)   
Aramova, Crona, and
De Negri define a notion of regularity based on the
degrees appearing in a free bigraded resolution of a finitely generated
module over a bigraded polynomial ring in \cite{ACD}.  (See also \cite{R}.)
We extend this notion to the more general case in which $M$ is a finitely generated multigraded module over the $\Z^k$-graded homogeneous coordinate ring of $\pnk$ by assigning to $M$ a \emph{resolution regularity vector}
$\underline{r}(M) \in \N^k$ (Definition \ref{def: res-reg})
which gives bounds on the degrees of the generators of the multigraded syzygy modules of $M.$ 

By contrast, Maclagan and
Smith \cite{MS} use the local cohomology definition of
regularity 
as their starting point in defining multigraded regularity for
toric varieties.  Since $\pnk$ is a toric variety, their definition
specializes to a version of multigraded regularity given
in terms of the vanishing of $H_{B}^i(M)_{\underline{p}}$,
the degree $\underline{p} \in \N^k$ part of the $i$th local
cohomology module of $M$, where $B$ is the irrelevant
ideal of the homogeneous coordinate ring of $\pnk.$
In this
context, the multigraded regularity, which we will denote by $\reg_B(M),$ 
is a subset of $\Z^k$.  It can be shown
that the degrees of the generators of $M$ lie outside the set $\bigcup_{1\leq i \leq k} (e_i + \reg_B(M)) ,$ and if $M$ is actually
$\N^k$-graded then the degrees of the generators must lie in $\N^k.$

If $k \geq 2$, the complement of  $\bigcup_{1\leq i \leq k} (e_i + \reg_B(M))$ in $\N^k$ may be unbounded; we thus lose a useful
 feature of regularity in the standard graded case.  However, when $k=2$,
Hoffman and Wang \cite{HW} introduced a notion of 
strong regularity (their definition of weak regularity essentially agrees
with \cite{MS}) which requires the vanishing of graded pieces of additional 
local cohomology modules and gives a bounded subset of
$\N^2$ 
that contains the degrees of the generators of $M.$  It would be interesting to develop a notion of 
strong regularity for other multigraded rings.  Since the completion of this paper, the authors,
together with Wang,
have developed a way of thinking about multigraded regularity via $\Z$-graded coarsenings of $\Z^r$-gradings.  (See \cite{SVW} for details.)

We now give an outline of the paper and describe our results.  In
\S 1 we briefly introduce multigraded
regularity as defined in \cite{MS}.
We also recall basic notions related to
fat point schemes in $\pnk$ and show that the degree of a fat point
scheme can be computed directly from the multiplicities of the points in Proposition \ref{prop: degree}.

In \S 2,  Proposition \ref{prop: res-reg} shows how
$\underline{r}(M)$ can be used to find a large subset of $\reg_B(M)$.
We also extend the work of \cite{R} to show how to use an almost
regular sequence to compute $\underline{r}(M)$ if $M$ is generated in degree $\underline{0}$ 
 in Theorem \ref{resregvec}.  We study the connections
between $\underline{r}(M)$ and the $\Z$-regularity of modules associated
to the factors of $\pnk$ in \S 3.

In \S 4 and \S 5, we study $\reg_B(R/I_Z)$ when
$R/I_Z$ is the coordinate ring of a fat point scheme $Z$ in $\pnk$.
Theorem \ref{thm: points res-vector} shows that
 $\underline{r}(R/I_Z) = (r_1,\ldots,r_k)$ where
$r_i = \reg(\pi_i(Z)) \subseteq \pr^{n_i}$.  Furthermore,
we show that $\underline{r}(R/I_Z) + \N^k \subseteq \reg_B(R/I_Z)$
in Proposition \ref{regbound}, which improves on bounds that follow from
Proposition \ref{prop: res-reg}.  
Moreover, if $Z \subseteq \pnk$ is also arithmetically Cohen-Macaulay,
we show that we have equality of sets:  $\underline{r}(R/I_Z)
+ \N^k = \reg_B(R/I_Z)$ (Theorem \ref{them: acm points}).  In \S 5, we restrict our attention
to fat points $Z = m_1P_1 + \cdots + m_sP_s$ in $\pr^1 \times \pr^1$ with support in generic position
and combine results of \cite{DS} and \cite{GuVT} to show that
\[
\{(i,j) \in \N^2 ~|~ (i,j) \geq (m_1-1,m_1-1) 
~\mbox{and}~ i+j \geq \max\{m-1,2m_1-2\}\} \subseteq \reg_B(Z).
\]
where $m = \sum m_i$ and $m_1 \geq m_2 \geq \cdots \geq m_s$ (Theorem 
\ref{them: p1xp1bound}).

At the end of the paper we have included an appendix containing some modified versions
of results found in \cite{MS}.

\noindent
{\it Acknowledgments.}  The authors would like to thank David Cox, T\`ai H\`a,
Diane Maclagan, Tim R\"omer, Greg Smith, Haohao Wang and an anonymous referee 
for their comments and suggestions, and the organizers of the COCOA VIII
conference for providing us with the opportunity to meet and begin this project.
The computer software packages CoCoA \cite{Co} and Macaulay 2 \cite{GS} were very helpful in computing examples throughout the project.
The first author would like to thank the Clare Boothe Luce Program for financial
 support.  The second author  also acknowledges the financial
support of NSERC.


\section{Setup}
\subsection{The homogeneous coordinate ring of $\pnk$}
Let {\bf k} be an algebraically closed field of characteristic zero. 
 Let $\N$ denote the natural numbers $0, 1, \ldots.$
The coordinate ring of $\pnk$ is the multigraded polynomial ring $R = {\bf k}[x_{1,0},\ldots,x_{1,n_1},
\ldots,x_{k,0},\ldots,x_{k,n_k}]$ where $\deg x_{i,j} = e_i$,
the $i$th standard basis vector of $\Z^k.$  Because $R$ is an $\N^k$-graded ring,  $R = \bigoplus_{\bi\in\N^k}
R_{\bi}$ and $R_{\bi}$ is a finite dimensional vector space over
{\bf k} with  a basis consisting of all  monomials of multidegree $\bi$.  Thus,
$\dim_{\bf k} R_{\bi} = \binom{n_1+i_1}{n_1}\binom{n_2+i_2}{n_2} \cdots\binom{n_k+i_k}{n_k},$
where $\bi = (i_1,\ldots,i_k)$.

Note that $R$ is the homogeneous coordinate 
ring of $\pnk$ viewed as a toric variety  of dimension $N := n_1 + \cdots + n_k.$  (See \cite{C} for a
comprehensive introduction to this point of view.)  The homogeneous coordinate ring of a toric variety 
is modeled after the homogeneous coordinate 
ring of $\pr^n.$  The space $\pnk$ is the quotient of $\A^{N+k} - V(B),$ where $B = \bigcap_{i = 1}^k \langle x_{i,j} \mid j = 0, \ldots, n_i \rangle$ is its square-free monomial ``irrelevant'' ideal. 
 Note that if $k = 1,$ then $B$ is just the
irrelevant maximal ideal of the coordinate ring of projective space.

The $\N^k$-homogeneous ideals of $R$ define subschemes of $\pnk.$ As in the standard
graded case, the notion of saturation plays an important role. 

\begin{definition}
Let $I \subseteq R$ be an $\N^k$-homogeneous ideal.  The saturation of $I$ with
respect to $B$ is 
$\sat(I) =  \{f \in R \mid f B^j \subseteq I, ~j \gg 0\}.$
\end{definition}

\noindent Two homogeneous ideals define the same subscheme of $\pnk$ if and only if
their saturations with respect to the irrelevant ideal are equal.  (See
Corollary 3.8 in \cite{C}.)

\subsection{Multigraded modules and regularity}
We shall work throughout with finitely generated
$\Z^k$-graded $R$-modules $M = \bigoplus_{\underline{t} \in \Z^k} M_{\underline{t}}.$  
Without loss of generality, we may restrict our
attention to $\N^k$-graded modules, since the $\underline{t} \in \Z^k$ with
$M_{\underline{t}} \neq 0$ 
must be contained in $\underline{p} + \N^k$ for some
$\underline{p} \in \Z^k$ if $M$ is finitely generated. Write 
$\underline{p} = \underline{p}^+ - \underline{p}^-$ where $\underline{p}^+, \underline{p}^- \in \N^k.$
Shifting degrees by $-\underline{p}^-$
yields a finitely generated $\N^k$-graded module.  

When $M$ is a finitely generated $\N^k$-graded $R$-module, it is useful to view $M$ as
both an $\N^1$-graded module and an $\N^k$-graded module.  We
introduce some notation and conventions for translating between the
$\N^k$ and $\N^1$ gradings 
of a module.  Let $\underline{a} = (a_1, \ldots, a_k) \in \N^k$, and let $\underline{1} = (1, \ldots, 1).$
If $m \in M$ has multidegree $\underline{a} \in \N^k,$ 
define its $\N^1$-degree to be $\underline{a}\cdot \underline{1}.$

We will use 
$\HH_{M}$ to denote the multigraded Hilbert function
$\HH_{M}(\bt) := \dim_{\bf k} M_{\bt}$, and $H_{M}$ to denote
the $\N^1$-graded Hilbert function $H_{M}(t) := \dim_{\bf k} M_t$.
Because $M_t = \bigoplus_{t_1+\cdots+t_k = t} M_{(t_1,\ldots,t_k)}$,
we have the identity
\[
H_{M}(t) = \sum_{t_1+\cdots+t_k=t} \HH_{M}(t_1,\ldots,t_k)
~\hspace{.5cm}\mbox{for all $t \in \N$.}
\]
If $I_Y$ is the $B-$saturated ideal defining 
a subscheme $Y \subseteq \pnk$, then we
sometimes write $\HH_{Y}$ (resp. $H_Y$) for $\HH_{R/I_Y}$ (resp. $H_{R/I_Y}$).

We use the notion of multigraded regularity developed in \cite{MS}.  
To discuss this notion of regularity, we require a
preliminary definition.  

\begin{definition}\label{defn: regions}
Let $i \in \Z$ and set
\[
\N^k[i] := \bigcup (\sgn(i) \underline{p}+\N^k) \subset \Z^k
\]
where the union is over all $\underline{p} \in \N^k$  whose coordinates sum to $|i|.$
(In the notation of \cite{MS}, \S 4, we have taken $\mathcal{C}$ 
to be the set of standard basis vectors of $\Z^k.$)
\end{definition}

Note that $\N^k[i]$ may not be contained in $\N^k.$  
The generality of Definition \ref{defn: regions} is necessary
because $\N^k[i]$ will be used to describe the degrees in which
certain local cohomology modules of $\N^k$-graded modules vanish, and these 
local cohomology modules may be nonzero in degrees with negative coordinates.

\begin{definition}[Definition 4.1 in \cite{MS}]\label{def: MS-reg}
Let $M$ be a finitely generated $\N^k$-graded $R$-module.
If $\underline{m} \in \Z^k$, we say that $M$ is $\underline{m}$-regular if
$H^i_B(M)_{\underline{p}} = 0$ for all $\underline{p} \in \underline{m}+\N^k[1-i]$ 
for all $i \ge 0.$ The multigraded regularity of $M,$ denoted $\reg_B(M),$ is the set 
of all $\underline{m}$ for which $M$ is $\underline{m}$-regular. 
\end{definition}

When $M = R/I_Y$, the $\N^k$-graded coordinate ring associated to a scheme
$Y \subseteq \pnk$, we shall write $\reg_B(Y)$ to denote $\reg_B(R/I_Y)$.

If $k=1$ and $M \neq 0,$ then $\reg_B(M)$ is a subset of $\N^1$, and there
exists some $r \in \N^1$ such that $\reg_B(M) = \{i ~|~ i \geq r\}$.
In this case, we will simply write $\reg(M) = \reg_B(M) = r$.  Note
that $\reg(M)$ is the standard Castelnuovo-Mumford regularity.
When $k = 2,$ Definition \ref{def: MS-reg} is essentially 
the same as 
the notion of \emph{weak regularity} (Definition 3.1 in \cite{HW})
of Hoffman and Wang.

\begin{remark}\label{lem: reg_B(R)}
As one might expect, $\reg_B(R) = \N^k.$ Indeed, it follows from Example 6.5
in \cite{MS} that it is enough to show the corresponding fact for a notion
of multigraded regularity for the sheaf $\mathcal{O}_{\pnk}.$   This can be done using
the K\"unneth formula generalizing the 
proof for the case $k = 2$ in Proposition 2.5 of \cite{HW}.  A topological approach is given in Proposition 6.10 of \cite{MS}.  For a related result see Proposition 4.3 in \cite{HW}.
\end{remark}

It will also be useful to have the following weaker condition of
\emph{multigraded regularity from level $\ell$}:

\begin{definition}[Definition 4.5 in \cite{MS}]\label{def: level l reg}
Given $\ell \in \N,$ the module $M$ is $\underline{m}$-regular from level $\ell$
if $H^i_B(M)_{\underline{p}} = 0$ for all $i\ge \ell$ and all $\underline{p} \in 
{\underline{m}}+\N^k[1-i].$  The set of all ${\underline{m}}$ such that $M$ is 
${\underline{m}}$-regular from level $\ell$ is denoted $\reg_B^{\ell}(M).$
\end{definition}

Note that $\reg_B^{\ell}(M) \supseteq \reg_B(M)$ for any finitely generated 
multigraded $R$-module $M.$  However,
even when $M = R,$ the inequality may be strict.

\begin{example}\label{ex: reg R}
Let $R$ be
the homogeneous coordinate ring of $\pr^2 \times \pr^2.$  We will show that $(-1,0) \in \reg_B^4(R).$  
By Definition \ref{def: level l reg} we only need to check vanishings of graded pieces of $H^i_B(R)$ for $i \ge 4$ which
is equivalent to checking vanishings of $H^i(\pr^2 \times  \pr^2, \mathcal{O}_{\pr^2 \times \pr^2}(a,b))$ for $i \ge 3.$  (See \S 6 in \cite{MS}.) 
We let $H^i(\mathcal{O}_{\pr^2 \times \pr^2}(a,b))$ denote  $H^i(\pr^2 \times \pr^2, \mathcal{O}_{\pr^2 \times \pr^2}(a,b))$ and
$H^i(\mathcal{O}_{\pr^2}(a))$ denote  $H^i(\pr^2, \mathcal{O}_{\pr^2}(a)).$
By the K\"unneth formula, $H^3(\mathcal{O}_{\pr^2 \times \pr^2}(a,b))$ is the direct sum 
\[ \bigoplus_{i+j = 3, \ i,j \ge 0} H^i(\mathcal{O}_{\pr^2}(a))\otimes H^j(\mathcal{O}_{\pr^2}(b)).\]
Since 
\[H^1(\mathcal{O}_{\pr^2}(d)) = H^3(\mathcal{O}_{\pr^2}(d)) = 0\] for all integers $d,$ each of the terms in the 
direct sum has a factor that is zero.  

Similarly, if we compute $H^4(\mathcal{O}_{\pr^2 \times \pr^2}(a,b))$ using the K\"unneth formula, the 
only possible nonzero contribution to the direct sum comes from 
$H^2(\mathcal{O}_{\pr^2}(a))\otimes H^2(\mathcal{O}_{\pr^2}(b))$, which is nonzero if and only if both $a,b \leq -3.$  
However, the vanishing conditions needed for $(-1,0)$ to be in $\reg_B^4(R)$ only require 
vanishing of $H^4(\mathcal{O}_{\pr^2 \times \pr^2}(a,b))$ for $(a,b) \ge (-5,0), (-4,-1), (-3,-2), (-2,-3), (-1, -4).$  

All of the cohomology 
groups $H^i(\mathcal{O}_{\pr^2 \times \pr^2}(a,b))$ vanish for $i \ge 5$ since
\[H^i(\mathcal{O}_{\pr^2 \times \pr^2}(a,b)) = \bigoplus_{j_1+j_2 = i} H^{j_1}(\mathcal{O}_{\pr^2}(a))\otimes H^{j_2}(\mathcal{O}_{\pr^2}(b))\] and $i \ge 5$ implies that at least one of
$j_1$ and $j_2$ is at least 3.
Therefore, $(-1,0) \in \reg_B^4(R) \supsetneq \reg_B(R).$
\end{example}

\subsection{Hilbert functions of points}
We recall some facts about points in multiprojective spaces.  If
$P \in \pnk$ is a point, then the ideal
$I_P \subseteq R$ associated 
to $P$ is the prime ideal
$I_P = \langle L_{1,1},\ldots,L_{1,n_1},\ldots,L_{k,1},\ldots,
L_{k,n_k}\rangle$ with $\deg L_{i,j} = e_i$. 
Let $X = \{P_1,\ldots,P_s\}$ be a set of distinct points in $\pnk$,
and let $m_1,\ldots,m_s$ be positive integers.  Set 
$\Iz = I_{P_1}^{m_1} \cap \cdots \cap I_{P_s}^{m_s}$
where $I_{P_i} \leftrightarrow P_i$.  Then $I_Z$ defines the scheme
of fat points $Z = m_1P_1 + \cdots + m_sP_s \subseteq \pnk$.

The degree of $Z$ is its length as a 0-dimensional subscheme of $\pnk.$

\begin{proposition}\label{prop: degree}
The degree of $Z = m_1P_1 + \cdots + m_sP_s$ is
\[ \sum_{i = 1}^s \binom{N+m_i-1}{m_i-1}.\] 
\end{proposition}

\begin{proof}
The ideal $I_Z$ is a $B$-saturated ideal defining a finite length subscheme 
of $\pnk.$  We will compute the degree of $Z$ by computing the lengths of
the stalks of the structure sheaf of $Z$ at each of the points $P_i.$ 

The stalk of $\mathcal{O}_{\pnk}$ at a point $P_i$ is a local ring isomorphic to 
$O = k[x_1, \ldots, x_N]_{\langle x_1, \ldots, x_N \rangle}$ where
the $x_i$ are indeterminates and  $\langle x_1, \ldots, x_N \rangle$
is a  maximal ideal.  The length of
$\mathcal{O}_{Z, P_i} = O/ \langle x_1, \ldots, x_N \rangle^{m_i}$ is 
\[\sum_{j = 0}^{m_i-1} {N +j-1\choose j},\] so the result follows
once we apply the identity $\sum_{k=0}^r \binom{n+k}{k} = \binom{n+r+1}{r}$.
\end{proof}

Short exact sequences constructed by taking a hyperplane section arise
frequently in proofs involving regularity in the standard graded case.  In the
multigraded generalization, we will employ the use of hypersurfaces  of each 
multidegree $e_i.$  
Algebraically, we need the following lemma, which generalizes the
reduced case of Lemma 3.3 in \cite{VT}.

\begin{lemma}\label{singlenzd}
If $I_Z$ is the defining ideal of $Z$, a set of fat points in
$\pnk,$ then, for each $i=1,\ldots,k$,  there exists an $L \in R_{e_i}$ that 
is a nonzerodivisor on $R/\Iz$. 
\end{lemma}

\begin{proof}
We will show only the case $i =1$.
Since the primary decomposition of $I_Z$ is $I_Z = I_{P_1}^{m_1} \cap \cdots \cap I_{P_s}^{m_s}$,
the set of zerodivisors of $R/I_Z$, consists of the set $\bigcup_{i=1}^s I_{P_i}.$  It will suffice to show ${\displaystyle \bigcup_{i=1}^s (I_{P_i})_{e_1} 
\subsetneq R_{e_1}}$.  It is clear that 
$(I_{P_i})_{e_1} \subsetneq R_{e_1}$ for each $i = 1,\ldots,s$.
Because the field ${\bf k}$ is infinite, the vector space $R_{e_1}$ 
cannot be expressed as a finite union of vector spaces, and hence,  
$\bigcup_{i=1}^s (I_{P_i})_{e_i} \subsetneq R_{e_1}$.
\end{proof}

Using Lemma \ref{singlenzd} we can describe rules governing the behavior of the
multigraded Hilbert function of a set of fat points.

\begin{proposition} \label{HFprop}
Let $Z$ be a set of fat points of $\pnk$ with Hilbert function $\HH_Z$.  Then
\begin{enumerate}
\item[$(i)$] for all $\bi \in \N^k$ and all $1 \leq j \leq k$, $\HH_Z(\bi) \leq \HH_Z(\bi + e_j)$.
\item[$(ii)$] if $\HH_Z(\bi) = \HH_Z(\bi+e_j)$, then $\HH_Z(\bi+e_j) = \HH_Z(\bi+2e_j)$.
\item[$(iii)$] $\HH_{Z}(\bi) \leq \deg(Z)$ for all $\bi \in \N^k$.
\end{enumerate}
\end{proposition}

\begin{proof}
To prove $(i)$ and $(ii)$ use the nonzerodivisors of Lemma \ref{singlenzd} 
to extend the proofs 
of Proposition 1.3 in \cite{GuVT} for fat points in $\popo$ to $\pnk$.
For $(iii)$, if $Z = m_1P_1+\cdots+m_sP_s$, then
$\deg (Z)$ is an upper bound on the number of linear conditions
imposed on the forms that pass through the points $P_1,\ldots,P_s$
with multiplicity at least $m_i$ at each point $P_i$.
\end{proof}

If $I_Z$ defines a set of fat points in $\pnk$, then 
the computation of $\reg_B(Z)$, as defined
by Definition \ref{def: MS-reg}, depends only upon knowing 
$\HH_Z$.  Indeed

\begin{theorem}[{Proposition 6.7 in \cite{MS}}]
\label{HF=reg} Let $Z$ be a set of fat points in $\pnk$.
Then $\bi \in \reg_B(Z)$ if and only
if $\HH_Z(\bi) = \deg(Z)$.
\end{theorem}

\begin{remark}\label{generic}
The set of reduced points $Z = P_1 + P_2 + \cdots + P_s$ 
is said to be in {\it  generic position}  
if $\HH_{Z}(\bi) = \min\{\dim_{\bf k} R_{\bi},s\}$ for all $\bi \in \N^k$.  Hence,
if $Z$ is in generic position, 
$\reg_B(Z) = \{\bi ~|~ \dim_{\bf k} R_{\bi} \geq s\}$.
\end{remark}


\section{Regularity and syzygies}
In the $\N^1$-graded case, the definition of regularity can be formulated in terms of the degrees that appear as 
generators in the minimal free graded 
 resolution of $M.$  In this section we discuss 
a multigraded version of this
definition extending the bigraded generalization that
was given in \cite{ACD} and studied further in \cite{R}.

We define a multigraded version of the notions of
$x$- and $y$-regularity from \cite{ACD}.

\begin{definition}\label{def: res-reg}
Let \[ r_{\ell} := \max \{a_{\ell}~|~\tor_i^R(M, {\bf k})_{(a_1, \ldots, a_{\ell}+i, \ldots, a_k)} \neq 0  \}\]
for some $i$ and for some $a_1, \ldots, a_{\ell-1}, a_{\ell+1}, \ldots, a_k.$  We
will call $\underline{r}(M) := (r_1, \ldots, r_k)$ the \emph{resolution regularity vector} of $M.$\end{definition}

Note that if $\underline{r}(M) = (r_1, \ldots, r_k)$ is the resolution regularity vector
 of a module $M,$ then the multidegrees appearing at the $i$th stage in the minimal
graded 
free resolution of $M$ have $\ell$th coordinate bounded above by $r_{\ell}+i.$  Indeed,  
\[R(-b_{1}, \ldots, -b_{k})_{(a_1, \ldots, a_{\ell}+i, \ldots, a_k)} \neq 0
\]
exactly when \[ (a_1-b_{1}, \ldots, a_{\ell}+i - b_{\ell}, \ldots, a_k-b_k)\]
has nonnegative coordinates, and 
\[ \tor_i^R(M, {\bf k})_{(a_1, \ldots, a_{\ell}+i, \ldots, a_k)} \neq 0\]
when $R(-a_{1}, \ldots, -a_{\ell-1}, -a_{\ell} - i, -a_{\ell+1}, \ldots, -a_{k})$
appears as a summand of the module at the $i$th stage in the resolution.

The resolution regularity vector of a module allows us to compute a
lower bound on the multigraded regularity of a module.

\begin{proposition}\label{prop: res-reg}
Let $M$ be a finitely generated $\N^k$-graded $R$-module with
resolution regularity vector $\underline{r}(M) = (r_1,\ldots,r_k)$.  
If 
$\underline{p} \in \underline{r}(M) + \N^k,$ then
\[\bigcup_{\underline{a} \in \N^k, 
\underline{a}\cdot\underline{1}= m-1} \underline{p} + m \cdot\underline{1}- \underline{a} +\N^k \subseteq \reg_B(M)\]
where $m = \min \{ N+1, \prdim M\}.$  Note that this set equals
\[\underline{p} + m \cdot \underline{1} + \N^k[-(m-1)]\] if $m > 0.$
\end{proposition}

\begin{proof}
Let $E.$ be the 
minimal free multigraded 
resolution of $M$ where $E_i = \bigoplus R(-{\underline{q}}_{ij})$ with 
$\underline{q}_{ij} \in \N^k$ and 
$\underline{q}_{ij} \leq \underline{p}+ i\cdot\underline{1}$ for $i =0,\ldots, \prdim M.$
Therefore, by Remark \ref{lem: reg_B(R)} we have the following bound 
on the multigraded regularity of $E_i,$ 
\[\reg_B(E_i) \supseteq \bigcap (\underline{q}_{ij} +\N^k )\supseteq \underline{p} + i \cdot \underline{1}+ \N^k.\]

If $\prdim M = 0,$ then it is immediate that $\underline{p}+\N^k \subseteq \reg_B(M).$  So assume $\prdim M > 0$.
Let $K_0$ be the syzygy module of $M$.  So we have a short exact sequence
\[0 \rightarrow K_0 \rightarrow E_0 \rightarrow M \rightarrow 0.\]
By Lemma \ref{appendix lemma} we then have 
\[\bigcap_{1 \leq j \leq k} (-e_j + \reg_B^1(K_0)) \cap \reg_B^0(E_0) \subseteq \reg_B^0(M) =
\reg_B(M).\]
Since  $\reg_B^1(K_0) \subseteq (-e_j + \reg_B^1(K_0))$ for
each $j$, we have 
\begin{eqnarray*}
\reg_B^1(K_0) \cap (\underline{p} + \N^k) &\subseteq &\reg_B^1(K_0) \cap \reg_B^0(E_0) \\
&\subseteq &\bigcap_{1 \leq j \leq k} (-e_j + \reg^1(K_0)) \cap \reg_B^0(E_0) \subseteq 
\reg_B(M).
\end{eqnarray*}
 
Applying the revised Corollary 7.3 of \cite{MS} (see Theorem \ref{appendix theorem}) implies that 
\begin{equation}\label{eqn: union1} \bigcup_{\phi:[m] \to [k]} 
\left( \bigcap_{ 1\leq i \leq m} -e_{\phi(2)} - \cdots - e_{\phi(i)} +
\reg^i(E_{i})\right) \subseteq \reg^1_B(K_0)\end{equation}
from which we deduce that
\begin{equation}\label{eqn: union} \bigcup_{\phi:[m] \to [k]} 
\left( \bigcap_{ 1\leq i \leq m} -e_{\phi(2)} - \cdots - e_{\phi(i)} +
\underline{p} + i \cdot \underline{1} + \N^k\right) \subseteq \reg^1_B(K_0).\end{equation}
Here we are using the fact that $\reg_B(E_i)$
is contained in $\reg_B^{\ell}(E_i)$ for all $\ell \ge 0.$

Suppose that $\phi:[m] \to [k].$  Consider 
\begin{equation} \label{eqn: intersection} 
\bigcap_{1 \leq i \leq m} -e_{\phi(2)} - \cdots - e_{\phi(i)} + \underline{p}+ i \cdot \underline{1}+ \N^k.
\end{equation}
The maximum value of the $j$th coordinate of 
$-e_{\phi(2)} - \cdots - e_{\phi(i)} + \underline{p}+ i \cdot \underline{1}$
over $i = 1, \ldots, m$ occurs when $i = m.$  Indeed, if the maximum value of the $j$th coordinate occurs for 
some $i < m,$ then consider 
$-e_{\phi(2)} - \cdots - e_{\phi(i)}-e_{\phi(i+1)} + \underline{p}+ (i+1)\cdot\underline{1}.$  
Since the difference 
\[ -e_{\phi(2)} - \cdots - e_{\phi(i)}-e_{\phi(i+1)} + \underline{p}+ (i+1)\cdot\underline{1} - (-e_{\phi(2)} - \cdots - e_{\phi(i)} + \underline{p}+ i\cdot\underline{1})\]
is
$-e_{\phi(i+1)} + \underline{1},$ we see that   
if $\phi(i+1) = j$, then the two vectors are equal in the $j$th coordinate.
Otherwise, the vector $-e_{\phi(2)} - \cdots - e_{\phi(i)}-e_{\phi(i+1)} + \underline{p}+ (i+1)\cdot\underline{1}$ has a bigger
$j$th coordinate.  So we see that the maximum value of each of the coordinates must occur when $i = m$ (and possibly earlier as well).

Therefore, the intersection in (\ref{eqn: intersection})
is equal to $\underline{p}+ m\cdot\underline{1}- \sum_{i = 2}^{m} e_{\phi(i)} + \N^k.$  As $\phi$ varies over all possible
functions from $[m]$ to $[k],$ the set of all vectors
$\sum_{i = 2}^{m} e_{\phi(i)}$ is just the set of all $\underline{a} \in \N^k$ such that
$\underline{a}\cdot\underline{1}= m-1.$  

(\ref{eqn: union}) is just 
\[\bigcup_{\underline{a} \in \N^k, \underline{a}\cdot\underline{1}= m-1} \underline{p} + m\cdot\underline{1}- \underline{a} +\N^k.\]
is also a subset of  $(\underline{p} + \N^k)$
we have 
\[\bigcup_{\underline{a} \in \N^k, \underline{a}\cdot\underline{1}= m-1} \underline{p} + m\cdot\underline{1}- \underline{a} +\N^k
\subseteq \reg_B^1(K_0) \cap (\underline{p}+\N^k) \subseteq \reg_B(M)\]
as desired.
\end{proof}

The $\N^1$-graded regularity of a multigraded module $M$ also gives the following rough bound on $\reg_B(M)$.

\begin{corollary}\label{cor: reg1-regk} 
Let $M$ be a finitely generated $\N^k$-graded $R$
module.  If $\reg(M) \leq r,$ then 
\[\reg_B(M) \supseteq 
 \bigcup_{\underline{a}\cdot\underline{1}= m-1}( r+m)\cdot\underline{1}-\underline{a}+\N^k\]
where $m = \min \{ N+1, \prdim M \}.$
\end{corollary} 
\begin{proof}   
Let $E.$ be the 
minimal free multigraded 
resolution of $M.$  Since $\reg(M) \leq r,$ we know that $\reg(E_i) \leq
r+i$ for all $i\ge 0.$   Since $E_i = 
\bigoplus R(-\underline{q}_{ij})$ with 
$\underline{q}_{ij} \in \N^k$, this means $\underline{q}_{ij}\cdot\underline{1}\leq r+i.$  Therefore, the resolution regularity vector $\underline{r}(M) 
\leq r \cdot \underline{1}.$  The result now follows from 
Proposition \ref{prop: res-reg}.
\end{proof}

\begin{remark}  Note that if $k = 1$ Proposition 
\ref{prop: res-reg} is equivalent to the statement that if an $\N^1$-graded
module is $p$-regular, then
it is also $q$-regular for all $q \in p + \N^1.$  

When $k >1,$ Proposition \ref{prop: res-reg}
may not give all of $\reg_B(M).$  Indeed, Theorem 7.2 of 
\cite{MS} will not give all of $\reg_B(M)$ even using more 
detailed information about multidegrees in a resolution than
given by $\underline{r}(M).$  (See 
Example 7.6 in \cite{MS}.) 
\end{remark}

\begin{remark}
Let $M$ be a finitely generated $\N^k$-graded $R$-module with resolution
regularity vector $\underline{r}(M)$.  Under extra
hypotheses on $M$, the bound of Proposition \ref{prop: res-reg} 
can be improved to $\reg_B(M) \supseteq \underline{r}(M) + \N^k$.
For example, in Proposition \ref{regbound} we will show that $\underline{r}(M) + \N^k
\subseteq \reg_B(M)$ if 
$M = R/I_Z$ is the coordinate ring of a set
of fat points in $\pnk.$

It seems, therefore,
natural to ask the following question:

\begin{question} \label{vector-reg}
Let $M$ be a finitely generated $\N^k$-graded $R$-module with resolution
regularity vector $\underline{r}(M)$.  What extra conditions on $M$ imply
$\reg_B(M) \supseteq \underline{r}(M) + \N^k?$  
\end{question}

(Since the submission of this paper, H\`a has given an example showing that this inclusion may not hold, as well as other
related results in \cite{ha}.)
\end{remark}

As we have observed, the resolution regularity vector of the $\N^k$-graded
$R$ module $M$
gives us partial information about $\reg_B(M)$.  We close this
section by describing how to compute the resolution regularity 
vector for some  classes of $M$.
This procedure is a natural extension of the bigraded case
as given by \cite{R}, which itself was a generalization of the
graded case \cite{AH}.

If $M$ is any finitely generated $\N^k$-graded $R$ module, then we shall
use $M_a^{[\ell]}$ to denote the $\N^{k-1}$-graded module
\[M_a^{[\ell]}:= \bigoplus_{(j_1,\ldots,j_{\ell-1},j_{\ell+1},\ldots,j_k)}
M_{(j_1,\ldots,j_{\ell-1},a,j_{\ell+1},\ldots,j_k)}.\]
Observe that $M_a^{[\ell]}$ is a ${\bf k}[x_{1,0},\ldots,
\hat{x}_{\ell,0},\ldots,\hat{x}_{\ell,n_{\ell}},\ldots,x_{k,n_k}]$-module,
where $\hat{~}$ means the element is omitted.

An element $x \in R_{e_{\ell}}$ is a {\it multigraded almost regular element} for $M$
if 
\[\left\langle 0 :_M x\right\rangle_a^{[\ell]} = 0 ~~\mbox{for $a \gg 0$.}
\]
A sequence $x_1,\ldots,x_t \in R_{e_{\ell}}$ is a 
{\it multigraded almost regular sequence}
if for $i=1,\ldots,t$, $x_i$ is a multigraded almost regular element for
$M/\langle x_1,\ldots,x_{i-1}\rangle M$.  A multigraded almost
regular element need not be almost regular in the usual sense, even 
for bigraded rings since we may have $\langle 0:_M  x \rangle_a^{[1]} =0$
for $a \ge a_0,$ but $\langle 0:_M x \rangle_{(a_0-1,j)} \neq 0$ for
infinitely many $j.$  (Note that in the single graded case, almost regular elements were
studied in \cite{T} under the name of filter regular elements.)

Now suppose that for each $\ell=1,\ldots,k$
we have a basis $y_{\ell,0}, \ldots, y_{\ell,n_{\ell}}$ of $R_{e_{\ell}}$ 
that forms a multigraded almost regular sequence for $M$.  Because ${\bf k}$ is infinite,
it is always possible to find such a basis; one can derive a proof of this
fact by adapting the proof of Lemma 2.1 of \cite{R} for the bigraded case to the multigraded case. 
Set
\[s_{\ell,j} := \max
\left\{a ~\left|~ \left\langle 0:_{M/\langle y_{\ell,0},\ldots,y_{\ell,j-1}\rangle M} y_{\ell,j} 
\right\rangle _a^{[\ell]} \neq 0\right\}\right.,\]
and $s_{\ell} := \max\{s_{\ell,0},\ldots,s_{\ell,n_{\ell}}\}$. 
Theorem 2.2 in \cite{R} then extends to the $\N^k$-graded case as follows:

\begin{theorem}\label{resregvec}
Let $M$ be a finitely generated multigraded $R$-module generated in degree $\underline{0},$
for for $\ell = 1,\ldots,k$, let $y_{\ell,0},\ldots,y_{\ell,n_{\ell}}$ be a basis
$R_{e_{\ell}}$ that forms a multigraded almost regular sequence for $M$. 
Then
$\underline{r}(M) = (s_1,\ldots,s_k)$.
\end{theorem}


\section{Resolution regularity and projections of varieties}

It is natural to ask if the $\N^1$-regularity of the projections of a subscheme $V$ of $\pnk$ 
onto the factors $\pr^{n_i}$ are related in a nice way to 
the coordinates appearing in the resolution regularity vector of $R/I_V.$
We show in Theorem \ref{thm: points res-vector} that if $V$ is a
set of fat points, then the $i$th coordinate in the resolution regularity vector of $R/I_V$ is
precisely the $\N^1$-regularity of the projection of $V$ to $\pr^{n_i}.$  
However, Example \ref{ex: projections} below shows that in general no such relationship
can hold for arbitrary subschemes of $\pnk.$  

\begin{example}\label{ex: projections}
Let $R = {\bf k}[x_0, x_1, x_2, y_0, y_1, y_2],$ and let 
\[I = \langle x_0, x_1 \rangle \cap \langle x_0-x_1, x_2\rangle \cap \langle y_0, y_1 \rangle \cap \langle y_0-y_1, y_2 \rangle\] be the defining ideal of a 
union of 4 planes in $\pr^2 \times \pr^2.$  The vector $\underline{r}(R/I)$ 
must be strictly positive in both coordinates since $I$ has a minimal generator
of bidegree $(2,2).$
However, the projection of the scheme onto either factor of $\pr^2$ is surjective.
Therefore, the regularity
of the projections of the scheme defined by $I$ is zero.
\end{example}

We consider some circumstances where the resolution regularity vector of a module $M$ is given by the regularities of modules associated to the factors of $\pnk.$  
We have the following proposition which generalizes Lemma 6.2 in \cite{R}.

\begin{proposition}\label{prop: product regularity}
Let $R_i = {\bf k}[x_{i,0}, \ldots, x_{i,n_i}]$ and let $M_i\neq 0$ be an $\N^1$-graded $R_i$-module.  The $i$th coordinate of the resolution regularity vector of $M_1 \otimes_{\bf k} \cdots \otimes_{\bf k} M_k$ is $\reg(M_i).$
\end{proposition}

\begin{proof}
The proof proceeds as in the case $k=2$ in \cite{R}. 
The point is that the tensor product (over {\bf k}) of minimal free graded resolutions
of the modules $M_i$ is the 
minimal free multigraded resolution of $M_1 \otimes_{\bf k} \cdots \otimes_{\bf k} M_k.$ 
We can read off the resolution regularity vector from the multidegrees appearing in this 
resolution.
\end{proof}

We have the following corollary.

\begin{corollary}\label{cor: product regularity}
Let $I_i$ be a proper homogeneous ideal in $R_i$, the coordinate 
ring of $\pr^{n_i}.$
Set $I := I_1 + \cdots + I_k \subset R$.  Then the $i$th coordinate of 
the resolution regularity vector of $R/I$ is $\reg(R_i/I_i).$
\end{corollary}
\begin{proof}
Let $M_i = R_i/I_i$ and apply Proposition \ref{prop: product regularity}.
\end{proof}

The $R$-modules $M$ that are products of modules over the factors of $\pnk$ have the 
property that 
they are $\underline{r}(M) \cdot \underline{1}$-regular as $\N^1$-graded modules.

\begin{corollary}\label{cor: tensor-reg}
Suppose that $M = M_1 \otimes_{\bf k} \cdots \otimes_{\bf k} M_k$ as in Proposition \ref{prop: product regularity} and  $\underline{r}(M) =(r_1, \ldots, r_k)$. Then $M$ 
is $\sum r_i$-regular as an $\N^1$-graded module. 
\end{corollary}

\begin{proof} Construct a resolution of $M$ by tensoring together 
minimal free graded resolutions of the $M_is.$  The free module at the $j$th
stage in the resolution is a direct sum of modules
$ F_{1,\ell_1}\otimes_{\bf k} \cdots \otimes_{\bf k} F_{k,\ell_k}$ 
where $\sum \ell_i= j$ and $F_{i,\ell_i}$ is the 
module at the $\ell_i$th stage in the 
minimal free graded resolution of $M_i$ over the ring $R_i$ defined as 
in Proposition \ref{prop: product regularity}. 
Since $F_{i,\ell_i}$ is generated by elements of degree $\leq r_i + \ell_i$,
 the total degree of any generator of $F_{i,\ell_1}\otimes_{\bf k} \cdots \otimes_{\bf k} F_{k,\ell_k}$
is $\leq \sum (r_i + \ell_i) = (\sum r_i) + j.$
\end{proof}

Corollary \ref{cor: tensor-reg} is not true for 
resolution regularity vectors of arbitrary modules.  For example,
let $R= {\bf k}[x_0, x_1, y_0, y_1]$ and let $I = \langle x_0y_1, x_1y_0 \rangle.$  (The saturation of $I$ with respect to $B$ is the defining ideal of two 
points in $\pr^1 \times \pr^1.$)  The resolution
of the ideal $I$ is given by the Koszul complex
\[0 \to R(-2,-2) \to R^2(-1,-1) \to \langle x_0y_1,x_1y_0 \rangle \to 0.\] 
So $\underline{r}(I) = (1,1)$, but the ideal $I$ is 
 not 2-regular as an $\N^1$-graded ideal.  

              
\section{Multigraded regularity for points}

Let $Z = m_1P_1 + \cdots + m_sP_s \subseteq \pnk$ be a scheme
of fat points, and let $Z_i = \pi_i(Z)$ denote the projection
of $Z$ into $\pr^{n_i}$ by the $i$th projection morphism
 $\pi_i:\pnk \rightarrow \pr^{n_i}$.  We show how
the resolution regularity vector of $R/I_Z$ is related to
$\reg(Z_i)$, the regularity of $Z_i$ as a subscheme
of $\pr^{n_i}$ for $i=1,\ldots,k$.
We then improve upon Proposition \ref{prop: res-reg}
and show $\underline{r}(R/I_Z) + \N^k \subseteq \reg_B(Z)$.
As a corollary, rough estimates of $\reg_B(Z)$ are obtained
for any set of fat points by
employing well known bounds for fat points in $\pr^n$.  We also
show that if $Z$ is ACM, $\reg_B(Z)$ is in fact determined by
$\reg(Z_i)$ for $i = 1,\ldots,k$.

\begin{lemma} \label{HFprojection}
Let $Z = m_1P_1 + \cdots + m_sP_s \subseteq \pnk.$ 
Let the $i$th coordinate of the $j$th point be $P_{ji}$ so that
the ideal $I_{P_j}$ defining the point $P_j$ is the sum 
of ideals $ I_{P_{j1}} + \cdots + I_{P_{jk}}$ where $I_{P_{ji}}$ defines
the $i$th coordinate of $P_j.$  Set $Z_i :=\pi_i(Z)$
for $i=1,\ldots,k$.  Then
\begin{enumerate}
\item[$(i)$]  $Z_i$ is the set of fat points in $\pr^{n_i}$
defined by the ideal $I_{Z_i} = \displaystyle \bigcap_{j = 1}^s I_{P_{ji}}^{m_j}.$
\item[$(ii)$] for all $t \in \N$, $H_{Z_i}(t) = \HH_Z(te_i)$.
\end{enumerate}
\end{lemma}

\begin{proof}
The proof of the reduced case found in
Proposition 3.2 in \cite{VT} can be adapted 
to the nonreduced case.
\end{proof}

\begin{theorem}\label{thm: points res-vector}
Suppose $Z \subseteq \pnk$ is a set of fat points.
Then
\[\underline{r}(R/I_Z) =  (r_1,\ldots,r_k)\] 
where $r_i = \reg(Z_i)$ for $i=1,\ldots,k$.
\end{theorem}

\begin{proof} 
We shall use Theorem \ref{resregvec} to compute the resolution regularity vector.
Because ${\bf k}$ is infinite, for each $\ell$ there exists a basis $y_{\ell,0},\ldots,
y_{\ell,n_{\ell}}$ for $R_{e_{\ell}}$ that is a multigraded almost regular
sequence.  Furthermore, by Lemma \ref{singlenzd} we can also assume that $y_{\ell,0}$ is
a nonzerodivisor on $R/I_Z$.

Since $y_{\ell,0}$ is a nonzerodivisor, $\langle 0:_{R/I_Z}y_{\ell,0}\rangle = 0$, which implies that $s_{\ell,0} \leq 0$.  
We now need to calculate $s_{\ell,i}$ for $i = 1,\ldots,
n_{\ell}$.

Because $y_{\ell,0}$ is 
a nonzerodivisor on $R/I_Z$, we have the short exact sequence
\begin{equation}\label{ses}
0 \rightarrow R/I_Z(-e_{\ell}) \stackrel{\times y_{\ell,0}}
{\rightarrow} R/I_Z \rightarrow R/\langle I_Z,y_{\ell,0}\rangle \rightarrow 0.
\end{equation}
Since $r_{\ell} = \reg(Z_{\ell})$, the sequence (\ref{ses}) and Lemma \ref{HFprojection}
give
\begin{eqnarray*}
\HH_{R/\langle I_Z,y_{\ell,0}\rangle}((r_{\ell}+1)e_{\ell}) &= &\HH_Z((r_{\ell}+1)e_{\ell})- \HH_Z(r_{\ell}e_{\ell})\\
& = & H_{Z_{\ell}}(r_{\ell}+1) - H_{Z_{\ell}}(r_{\ell}) = \deg Z_{\ell} - \deg Z_{\ell} = 0.
\end{eqnarray*}
Thus $\langle I_Z, y_{\ell,0}\rangle_{ae_{\ell}} = R_{ae_{\ell}}$ if $a \ge r_{\ell} +1$.  Hence
for any $\bj \geq (r_{\ell}+1)e_{\ell}$,
$R_{\bj} = \langle I_Z,y_{\ell,0}\rangle_{\bj} \subseteq \langle I_Z,y_{\ell,0},\ldots,y_{\ell,i-1}\rangle_{\bj}$.

Since $\langle 0 :_{R/\langle I_Z,y_{\ell,0},\ldots,y_{\ell,i-1}\rangle} y_{\ell,i}\rangle$ is an ideal
of $R/\langle I_Z,y_{\ell,0},\ldots,y_{\ell,i-1}\rangle$, and because 
$R/\langle I_Z,y_{\ell,0},\ldots,y_{\ell,i-1}\rangle_{\bj} = 0$, if $\bj \geq
(r_{\ell}+1)e_{\ell}$, 
\begin{equation}\label{ideal}
\left\langle 0 :_{R/\langle I_Z,y_{\ell,0},\ldots,y_{\ell,i-1}\rangle} y_{\ell,i} \right\rangle_a^{[\ell]} 
= 0 ~~\mbox{if $a \geq r_{\ell}+1$.}
\end{equation}
Thus, from (\ref{ideal}) we have $s_{\ell,j} \leq r_{\ell}$
for each $\ell$ and each $j=1,\ldots,n_{\ell}$.  Since $s_{\ell,0} \leq 0$,
it suffices to show that $s_{\ell,1} = r_{\ell}$
since this gives $s_{\ell} = \max\{s_{\ell,0},\ldots,s_{\ell,n_{\ell}}\} = s_{\ell,1} = r_{\ell}$.
The short exact sequence (\ref{ses}) also implies that
\[
\HH_{R/\langle I_Z,y_{\ell,0}\rangle}(r_{\ell}e_{\ell}) = H_{Z_{\ell}}(r_{\ell}) - H_{Z_{\ell}}(r_{\ell}-1) > 0\]
because $H_{Z_{\ell}}(r_{\ell}-1) < H_{Z_{\ell}}(r_{\ell}) = \deg Z_{\ell}$.
So there
exists $0 \neq \overline{F} \in (R/\langle I_Z,y_{\ell,0}\rangle)_{r_{\ell}e_{\ell}}$.  
Because
$\deg \overline{F}\overline{y}_{\ell,1} = (r_{\ell}+1)e_{\ell}$, and 
$(R/\langle I_Z,y_{\ell,0}\rangle)_{(r_{\ell}+1)e_{\ell}}
=0$, we must have 
$\overline{F} \in \langle 0 :_{R/\langle I_Z,y_{\ell,0}\rangle} y_{\ell,1}\rangle$.
So, $0 \neq \overline{F} \in \langle 0 :_{R/\langle I_Z,y_{\ell,0}\rangle} y_{\ell,1}\rangle_{r_{\ell}}^{[\ell]}$,
thus implying $s_{\ell,1} = r_{\ell}$.
\end{proof}

The previous result, combined with Proposition \ref{prop: res-reg},
gives us a crude bound on $\reg_B(Z)$.  However, 
we can improve upon this bound.

\begin{lemma} \label{HFfatpoint}
Let $P \in \pnk$ be a point with
defining ideal $I_P \subseteq R$, and $m \in \N^+$.  Then
$\reg_B(R/I_P^m) = (m-1,\ldots,m-1) + \N^k$.
\end{lemma}

\begin{proof} After a change of coordinates, we can assume
$P = [1:0:\cdots:0] \times \cdots \times [1:0\cdots:0]$.  So
$I_P^m = \langle x_{1,1},\ldots,x_{1,n_1},\ldots,x_{k,1},\ldots,x_{k,n_k}\rangle ^m.$
Since $I_P^m$ is a monomial ideal, $\HH_{R/I_P^m}(\bi)$
equals the number of monomials of degree $\bi$ in $R$ not
in $I_P^m$.  

A monomial $\prod x_{j,\ell}^{a_{j,\ell}}
\not\in (I_P^m)_{\bi}$ if and only
if  
$a_{1,1} + \cdots + a_{k,n_k} \leq m-1$ and 
$a_{j,1} + \cdots + a_{j,n_j} \leq i_j$ for each $j=1,\ldots,k$.
The result now follows since 
\begin{eqnarray*}
\#\left\{(a_{1,1},\ldots,a_{k,n_k}) \in \N^N ~\left|~ 
a_{1,1} + \cdots + a_{k,n_k} \leq m-1, \ \forall j \ a_{j,1} + \cdots + a_{j,n_j} \leq i_j 
\right\}\right.\end{eqnarray*}
is equal to
\begin{eqnarray*}
\#\left\{(a_{1,1},\ldots,a_{k,n_k}) \in \N^N ~\left|~ 
a_{1,1} + \cdots + a_{k,n_k} \leq m-1
\right\}\right.= \binom{m-1+N}{m-1} = \deg(mP)
\end{eqnarray*}
if and only if $\bi = (i_1,\ldots,i_k) \geq (m-1,\ldots,m-1)$.
\end{proof}

\begin{proposition}\label{regbound}
Suppose $Z \subseteq \pnk$ is a set of fat points. Then
\[(r_1,\ldots,r_k) + \N^k \subseteq \reg_B(Z)\]
where $r_i = \reg(Z_i)$ for $i = 1,\ldots,k$. 
\end{proposition}

\begin{proof} The proof is by induction on $s$, the
number of points in the support.  If $s = 1$, then the result follows from
Lemma \ref{HFfatpoint}.   

So, suppose $s > 1$ and let 
$X = \{P_1,\ldots,P_s\}$ be the support of $Z$.  We can
find an $i \in \{1,\ldots,k\}$ such that $1 < |\pi_i(X)|$, i.e.,
there exists an $i$ where the projection of $X$ onto
its $i$th coordinates consists of two or more points.  Fix a
$\tilde{P} \in \pi_i(X)$.  We can then write $Z = Y_1 \cup Y_2$
where $Y_1 = \{m_jP_j \in Z ~|~ \pi_i(P_j) = \tilde{P}\}$ and
$Y_2 = \{m_jP_j \in Z ~|~ \pi_i(P_j) \neq \tilde{P}\}$.  By
our choice of $i$, $Y_1$ and $Y_2$ are nonempty, $Y_1 \cap Y_2 
= \emptyset$, and $\pi_i(Y_1) \cap \pi_i(Y_2) = \emptyset$.

Let $I_{Y_1}$, resp., $I_{Y_2}$, denote the defining ideal associated
to $Y_1$, resp., $Y_2$.  Consider the short exact
sequence
\[
0 \rightarrow R/\langle I_{Y_1} \cap I_{Y_2}\rangle  \rightarrow R/I_{Y_1} \oplus
R/I_{Y_2} \rightarrow R/\langle I_{Y_1} + I_{Y_2}\rangle  \rightarrow 0.
\]
Since $I_Z = I_{Y_1} \cap I_{Y_2}$, this exact sequence gives rise
to the identity
\begin{equation} \label{eq: id}\HH_{Z}(\bt) = \HH_{Y_1}(\bt) + \HH_{Y_2}(\bt) - 
\HH_{R/\langle I_{Y_1}+I_{Y_2}\rangle }(\bt)~\mbox{for all $\bt \in \N^k$.}
\end{equation}
Set $Y_{j,1} := \pi_j(Y_1)$ and $Y_{j,2} :=\pi_j(Y_2)$ for
each $j = 1,\ldots,k$.  Since $Y_{j,1} \subseteq Z_j$ and 
$Y_{j,2} \subseteq Z_j$, we have $\reg(Y_{j,1}) \leq r_j$
and $\reg(Y_{j,2}) \leq r_j$.  By induction and the above
identity we therefore have
\[\HH_{Z}(r_1,\ldots,r_k) = \deg Y_1 + \deg Y_2 - 
\HH_{R/\langle I_{Y_1}+I_{Y_2}\rangle }(r_1,\ldots,r_k).\]
Since $\deg Z = \deg Y_1 + \deg Y_2$, 
it suffices to show $\HH_{R/\langle I_{Y_1} + I_{Y_2}\rangle }(r_1,\ldots,r_k) = 0$.

By Lemma \ref{HFprojection}, $\HH_{Z}(r_ie_i) = H_{Z_i}(r_i) = \deg Z_i, $
$\HH_{Y_1}(r_ie_i) = H_{Y_{i,1}}(r_i) = \deg Y_{i,1},$
and $\HH_{Y_2}(r_ie_i) = H_{Y_{i,2}}(r_i) = \deg Y_{i,2}$.
But because $Y_{i,1} \cap Y_{i,2} = \emptyset$ by our choice of $i$,
$\deg Z_i = \deg Y_{i,1} + \deg Y_{i,2}$.  Substituting
into (\ref{eq: id}) with $\underline{t} = r_ie_i$ then gives 
$\HH_{R/\langle I_{Y_1}+I_{Y_2}\rangle }(r_ie_i) = 0$, or equivalently,
$R_{r_ie_i} = \langle I_{Y_1} + I_{Y_2}\rangle _{r_ie_i}$.  It now follows
that $R_{(r_1,\ldots,r_k)} = \langle I_{Y_1}+I_{Y_2}\rangle _{(r_1,\ldots,r_k)}$
which gives $\HH_{R/\langle I_{Y_1} + I_{Y_2}\rangle }(r_1,\ldots,r_k) = 0$.
\end{proof}

Using well known bounds for fat points in $\pr^n$ thus gives us:

\begin{corollary}
Let $Z = m_1P_1+\cdots+m_sP_s \subseteq \pnk$ with $m_1 \geq \cdots \geq m_s$.
\begin{enumerate}
\item[$(i)$] Set $m=m_1+m_2+\cdots+m_s -1$.
Then $(m,\ldots,m) + \N^k \subseteq \reg_B(Z).$
\item[$(ii)$] Suppose $X = \{P_1,\ldots,P_s\}$
is in generic position.  For $i=1,\ldots,k$ set
\[\ell_i = \max\left\{m_1+m_2+1, \left\lceil\frac{(\sum_{i=1}^s m_i)+n_i-2}{n_i}\right\rceil
\right\}.\]
Then $(\ell_1,\ldots,\ell_k) + \N^k \subseteq \reg_B(Z).$
\end{enumerate}
\end{corollary}

\begin{proof}
It follows from Davis and Geramita \cite{DG}
that $r_i = \reg(Z_i) \leq m$ for each $i$.  
So $(m,\ldots,m) + \N^k \subseteq (r_1,\ldots,r_k) + \N^k$, and hence
$(i)$ follows.

For $(ii)$, because $X$ is in generic position, the support of $Z_i$ is in 
generic position in $\pr^{n_i}$.   In \cite{CTV} it was shown that 
$r_i = \reg(Z_i) \leq \ell_i$ for each $i$.
\end{proof}

Recall that a scheme $Y \subseteq \pnk$ is {\it arithmetically Cohen-Macaulay} 
(ACM) if $\operatorname{depth} R/I_Y = $K-$\dim R/I_Y$.  For any
collection of fat points $Z \subseteq \pnk$ we always have K-$\dim R/I_Z = k$,
the number of projective spaces.  However, for each $\ell \in \{1,\ldots,k\}$
there exist  sets of fat points (in fact, reduced points) $X_{\ell} \subseteq \pnk$
with $\operatorname{depth} R/I_{X_{\ell}} = 
\ell$.
See \cite{VT1} for more details.

A scheme of fat points, therefore, may or may not be ACM.
When $Z$ is ACM, $\reg_B(Z)$ depends only upon knowing $\reg(Z_i)$ 
for $i = 1,\ldots,k$.  

\begin{lemma}\label{nzd}
Let Z be an ACM set of fat points in $\pnk$.  Then there exist elements
$L_i \in R_{e_i}$ such that $L_1, \ldots, L_k$ is a regular sequence on 
$R/I_Z.$
\end{lemma}

\begin{proof}
The nontrivial part of the statement is the existence of a regular 
sequence whose elements have the specified multidegrees.  The proof
given for the reduced case (see Proposition 3.2 in \cite{VT1}) can be adapted
to the nonreduced case.
\end{proof}

\begin{theorem}\label{them: acm points}
Let $Z \subseteq \pnk$ be a set of fat points.  If $Z$ is ACM, 
then
\[(r_1,\ldots,r_k) + \N^k = \reg_B(Z)\]
where $r_i = \reg(Z_i)$ for $i = 1,\ldots,k$.
\end{theorem}
\begin{proof}
Let $L_1,\ldots,L_k$ be the regular sequence from Lemma \ref{nzd},
and set $J = \langle I_Z,L_1,\ldots,L_k\rangle $.
We require the following claims.

\noindent
{\it Claim 1. } If $\bj \not\leq (r_1,\ldots,r_k),$ then $\HH_{R/J}(\bj) = 0$. 

Since $\bj \not\leq (r_1,\ldots,r_k)$ there exists $1 \leq \ell \leq k$ 
such that $j_{\ell} > r_{\ell}$.  Using the exact sequence (\ref{ses})
of Theorem \ref{thm: points res-vector}, the claim
follows if we replace $y_{\ell,0}$ with $L_{\ell}$.

\noindent
{\it Claim 2.}  For $i=1,\ldots,k$, $\HH_{R/J}(r_ie_i) > 0$.

By degree considerations, $\HH_{R/J}(r_ie_i) = \HH_{R/\langle I_Z,L_i\rangle }(r_ie_i)$ for
each $i$.  Employing the short exact sequence
\[
0 \rightarrow R/I_Z(-e_i) \stackrel{\times L_i}{\rightarrow}
R/I_Z \rightarrow R/\langle I_Z,L_i\rangle  \rightarrow 0
\]
to calculate $\HH_{R/\langle I_Z,L_i\rangle }(r_ie_i)$ gives
$\HH_{R/\langle I_Z,L_i\rangle }(r_ie_i) = \HH_Z(r_ie_i) - \HH_Z((r_i-1)e_i) = H_{Z_i}(r_i) 
- H_{Z_i}(r_i-1)$.  The claim now follows since 
$H_{Z_i}(r_i-1) < \deg Z_i = H_{Z_i}(r_i)$.

We complete the proof.   Since $L_1,\ldots,L_k$ is a regular sequence,
we have the following short exact sequences
\[
0 \rightarrow R/\langle I_Z,L_1,\ldots,L_{i-1}\rangle (-e_i) \stackrel{\times L_i}
{\rightarrow} R/\langle I_Z,L_1,\ldots,L_{i-1}\rangle  \rightarrow R/\langle I_Z,L_1,\ldots,L_i\rangle  
\rightarrow 0
\]
for $i=1,\ldots,k$.   It then follows that
\[ 
\HH_{Z}(\bj) = \sum_{\underline{0} \leq \bi \leq \bj} \HH_{R/J}(\bi) ~~\mbox{for all $\bj \in \N^k$.}
\]
Now suppose that 
$\bj \not\in (r_1,\ldots,r_k) + \N^k$.  So $j_{\ell} < r_{\ell}$ for some $\ell$.
Set $j'_i = \min\{j_i,r_i\}$ and let $\bj' = (j'_1,\ldots,j'_k)$.  Note
that $\bj' \leq (r_1,\ldots,r_k)$ and $j'_{\ell} = j_{\ell} < r_{\ell}$.  By Claim 1 and the 
above identity
\[\HH_{Z}(\bj) = \sum_{\underline{0} \leq \bi \leq \bj} \HH_{R/J}(\bi) 
= \sum_{\underline{0} \leq \bi \leq \bj'} \HH_{R/J}(\bi).
\]
Then, by Claim 2,
\begin{eqnarray*}
\HH_{Z}(\bj) & = &\sum_{\underline{0} \leq \bi \leq \bj'} \HH_{R/J}(\bi) 
< \sum_{\underline{0} \leq \bi \leq \bj'} \HH_{R/J}(\bi) + H_{R/J}(r_{\ell}e_{\ell}) \\
&\leq& \sum_{\underline{0} \leq \bi \leq (r_1,\ldots,r_k)} \HH_{R/J}(\bi)
= \HH_Z(r_1,\ldots,r_k) = \deg Z.
\end{eqnarray*}
So, if $\bj \not\in (r_1,\ldots,r_k) + \N^k$, then $\bj \not\in \reg_B(Z)$.  This
fact, coupled with Theorem \ref{regbound}, gives the desired result.
\end{proof}

\begin{remark} The converse of Theorem \ref{them: acm points} is false
because there exist fat point schemes $Z$ such that 
$\reg_B(Z) = \underline{r}(R/I_Z) + \N^k$,
but $Z$ is not ACM.  For example, 
set $P_{ij} := [1:i] \times [1:j] \in \popo$, and consider
$Z = \{P_{11},P_{12},P_{13},P_{21},P_{22},P_{31},P_{33}\}$.
Then $\HH_Z$ is 
\[
\HH_{Z} = 
\bmatrix
1 & 2 & 3 & 3 & \ldots \\
2 & 4 & 6 & 6 & \ldots \\
3 & 6 & 7 & 7 & \ldots \\
3 & 6 & 7 & 7 & \ldots \\
\vdots&\vdots&\vdots&\vdots&\ddots
\endbmatrix.
\]
It follows that $\reg_B(Z) = (2,2) + \N^2$.  However,
 the resolution of $R/I_Z$ has length 3, so by  the Auslander-Buchsbaum Theorem,
depth $R/I_Z  = \operatorname{depth} R - \operatorname{pd} R/I_Z = 4 -3 =1 < 
2 = \operatorname{K-}\dim R/I_Z$.  Alternatively,
$R/I_{Z}$ is not ACM since
the first difference function of $\HH_{Z}$ is not the Hilbert function
of an $\N^2$-graded Artinian quotient of ${\bf k}[x_1,y_1]$. (see Theorem 4.8 in \cite{VT1}).
\end{remark}


\section{A bound for  fat points in $\popo$ with support
in generic position}

Let $Z = m_1P_1 + \cdots +m_sP_s$ be a set of fat points in $\popo$,
and furthermore, suppose that $X = \{P_1,\ldots,P_s\}$, the support of $Z$,
 is in generic
position, i.e.,
$\HH_X(i,j) = \min\{\dim_{\bf k} R_{(i,j)},s\}$ 
for all
$(i,j) \in \N^2$.  Using Proposition \ref{regbound}, we can obtain
the bound $(m-1,m-1) + \N^k \subseteq \reg_B(Z)$ where $m = \sum m_i$.
However, under these extra hypotheses, we can give a much
stronger bound.  

\begin{theorem}\label{them: p1xp1bound}
Let $Z = m_1P_1 + \cdots + m_sP_s \subseteq \popo$ be a set
of fat points whose support is in generic position.  
Assume $m_1 \geq m_2 \geq \cdots \geq m_s$, and set 
$m = m_1 + m_2 + \cdots +m_s$.  Then
\[
\{(i,j) \in \N^2 ~|~ (i,j) \geq (m_1-1,m_1-1) 
~\mbox{and}~ i+j \geq \max\{m-1,2m_1-2\}\} \subseteq \reg_B(Z).
\]
\end{theorem}

We require a series of lemmas.

\begin{lemma}\label{eventualgrowth}
Assume $Z \subseteq \popo$ is as in Theorem \ref{them: p1xp1bound}.
For $j=0,\ldots,m_1-1$, set 
\[c_j := \sum_{i=1}^s [m_i + (m_i-1)_+ + \cdots + (m_i-j)_+]
~\mbox{where $(n)_+ := \max\{0,n\}.$} \]
Then $c_{m_1-1} = \deg Z$.  As well, if we write $\HH_Z$ as
an infinite matrix, $\HH_Z$ has the following eventual behavior:
\[
\HH_{Z} = 
\begin{array}{ll}
 & {\scriptstyle \hspace{.4cm}0 \hspace{2.8cm} m_1-1\hspace{1.8cm} m-1} \\
\begin{array}{r}
{\scriptstyle 0} \\
\\
\\
\\
{\scriptstyle m_1-1}\\
\\
\\
{\scriptstyle m-1}\\
\\
\\
\end{array} 
&
\bmatrix
&  &   & &  &  & c_0 & c_0 & \cdots \\
&  &   & &  &  & c_1 & c_1 &\cdots \\
&  &   & &  &  & \vdots& &\\ 
&  & *  & &  &  & c_{m_1-2} & c_{m_1-2}& \cdots \\
&  &   & &  &  & \deg Z & \deg Z &\cdots \\
&  &   & &  &  & \vdots & &\\
c_0&c_1&\cdots &c_{m_1-2}&\deg Z & \cdots & \deg Z & \deg Z &\cdots \\
c_0&c_1&       &c_{m_1-2}&\deg Z &  & \deg Z & \deg Z & \cdots \\
\vdots&\vdots  &   &\vdots &\vdots  &  & \vdots&  \vdots&\ddots
\endbmatrix
\end{array}.
\]
\end{lemma}
\begin{proof}
Because $m_1\geq \cdots \geq m_s$ we have 
\[c_{m_1-1} = \sum_{i=1}^s [1+2+\cdots+m_i] = \sum_{i=1}^s \binom{m_i+1}{s}
= \deg Z.\]
The eventual behavior of $\HH_Z$ can be obtained
from Theorem 3.2 in \cite{GuVT}.
\end{proof}

\begin{lemma}\label{N^1prop}
Let $R/I_Z$ be the coordinate ring of a set of fat points that
satisfies
 the conditions of Theorem \ref{them: p1xp1bound}.
\begin{enumerate}
\item[$(i)$] As an $\N^1$-graded ring, $\reg(R/I_Z) \leq m-1$.
\item[$(ii)$] As an $\N^1$-graded ring, the Hilbert polynomial
of $R/I_Z$ is
\[HP_{R/{I_Z}}(t) = \sum_{i=1}^s \left[\binom{m_i+1}{2}t + \binom{m_i+1}{2}
\left(\frac{-2m_i + 5}{3}\right)\right].\]
\end{enumerate}
\end{lemma}
\begin{proof}
$(i)$ By using Theorem 4.4 in \cite{DS}, we obtain the bound
$\reg(I_Z) \leq  m.$ 
The result now follows since $\reg(R/I_Z) = \reg(I_Z) -1$.

$(ii)$ If $I_{P_i}$ is the defining ideal of a point $P_i$
in the support, then as an $\N^1$-homogeneous ideal of 
$R= {\bf k}[x_0,x_1,y_0,y_1]$, $I_{P_i}$ defines a line in $\pr^3$.  Since
the points in the support are in generic position, the lines that they
correspond to in $\pr^3$ must all be skew.

The structure sheaf $\mathcal{O}_Z = \bigoplus_{i = 1}^s \mathcal{O}_{m_iP_i}$, and hence, the
Hilbert polynomial of $\mathcal{O}_Z$ is just the sum of the Hilbert polynomials of $\mathcal{O}_{m_iP_i}.$  
The ideal $I_{P_i}^{m_i}$ is a power of complete intersection. The
resolution of $I_{P_i}^{m_i}$ is thus given by the Eagon-Northcott resolution.
Furthermore, since the generators of $I_{P_i}^{m_i}$ all have the
same degree, the Eagon-Northcott resolution produces
the following minimal graded free resolution of $R/I_{P_i}^{m_i}$:
\[0 \longrightarrow R^{m_i}(-(m_i+1)) \longrightarrow R^{m_i+1}(-m_i) 
\longrightarrow R \longrightarrow R/I^{m_i}_{P_i} 
\longrightarrow 0.\] 
We can compute the Hilbert polynomial of $R/I_{P_i}^{m_i}$ from
this resolution.  Since the Hilbert polynomial of $R/I_Z$ and its 
sheafification agree, we are done. 
\end{proof}

\begin{proof}{(of Theorem \ref{them: p1xp1bound})}
Set $\ell = \max\{m-1,2m_1-2\}$.
It suffices to show
that $\HH_Z(i,j) = \deg Z$ if $(i,j) \geq (m_1-1,m_1-1)$ and $i+j =\ell$.  
The conclusion then follows from
Proposition \ref{HFprop}.

By Lemma \ref{N^1prop}, $\reg(Z) \leq m-1$.  Thus $H_Z$,
the Hilbert function of $Z$ as a graded ring, agrees with
$HP_{R/I_Z}$ for all $t \geq m-1$.  In particular, since $\ell \geq m-1$,
\begin{equation}\label{HPrewrite}
H_Z(\ell) = \sum_{i=1}^s \left[\binom{m_i+1}{2}\ell + \binom{m_i+1}{2}
\left(\frac{-2m_i + 5}{3}\right)\right].
\end{equation}
Now $H_Z(\ell) = \sum_{i+j =\ell} \HH_Z(i,j)$.  If $i + j = \ell$, then
there are three cases:
\begin{enumerate}
\item[1.]  If $j \leq m_1 -2$, then by Lemma \ref{eventualgrowth}
and Proposition \ref{HFprop} (i), we have $\HH_Z(i,j) \leq c_j$.
\item[2.]  If $i \leq m_1 -2$, then we have $\HH_Z(i,j) \leq c_i$.
\item[3.]  If $j \geq m_1-1$ and $i \geq m_1 -1$, then $\HH_Z(i,j) \leq \deg Z$.
\end{enumerate}
Since there are $\ell+1$ pairs $(i,j) \in \N^2$ such that $i+j = \ell$,
and $m_1-1$ pairs fall into the first case, and $m_1-1$ are 
in the second case, we must have $\ell+1 - 2(m_1-1) = \ell+1 - 2m_1 +2
> 0$ pairs in the third case because $\ell \geq 2m_1-2$.

We thus have
\[H_Z(\ell) = \sum_{i+j = \ell} \HH_Z(i,j) \leq 2(c_0+\cdots+c_{m_1-2}) + 
(\ell+1 -2m_1+2)\deg Z.\]

\noindent
{\it Claim.} ${\displaystyle c_0 + \cdots + c_{m_1-2} = 
\sum_{i=1}^s \left[ \binom{m_i+1}{2}\left(\frac{2m_i-2}{3}\right) 
+ \binom{m_i+1}{2}(m_1-m_i)\right]}$.

\noindent
{\it Proof of Claim.}
\begin{eqnarray*}
\sum_{j=0}^{m_1-1} c_{j}& =& \sum_{j=0}^{m_1-1} \sum_{i=1}^s
[m_i + (m_i-1)_+ + \cdots + (m_i - j)_+] \\
&=& \sum_{i=1}^s \left[m_i^2 + (m_i-1)^2 + \cdots + 2^2 + 1^2
+ (m_1-m_i)\binom{m_i+1}{2} \right] \\
& = & \sum_{i=1}^s \left[\frac{m_i(m_i+1)(2m_i+1)}{6} + 
(m_1-m_i)\binom{m_i+1}{2} \right] \\
& = & \sum_{i=1}^s \left[\binom{m_i+1}{2}\left(\frac{2m_i+1}{3}\right) + 
(m_1-m_i)\binom{m_i+1}{2} \right]. 
\end{eqnarray*}
By subtracting $c_{m_1-1}
= \sum_{i=1}^s \binom{m_i+1}{2}$ from both sides of the
above expression, we arrive at the claimed result.

We complete the proof:
\begin{eqnarray*}
H_Z(\ell) & = & \sum_{i+j=\ell} \HH_Z(i,j) \\
&\leq & 2 \sum_{i=1}^s \left[\binom{m_i+1}{2}\left(\frac{2m_i-2}{3}\right) 
+ \binom{m_i+1}{2}(m_1-m_i)\right]  \\
&&+(\ell+1-2m_1+2)\sum_{i=1}^s \binom{m_i+1}{2}\\
& = & \sum_{i=1}^s \left[\binom{m_i+1}{2}\ell + \binom{m_i+1}{2}
\left(\frac{4m_i-4}{3} - 2m_i + 3\right)\right] \\
& = & \sum_{i=1}^s \left[\binom{m_i+1}{2}\ell + \binom{m_i+1}{2}
\left(\frac{-2m_i + 5}{3}\right)\right] = H_Z(\ell).
\end{eqnarray*}
If we view $\HH_Z$ as an infinite matrix, then because
$\HH_Z$ strictly increases along each row and column
until it reaches its eventual growth value as given in Lemma 
\ref{eventualgrowth},
we must have $\HH_Z(i,j)$ equal to this eventual growth value
for all $(i,j)$ with  $i+j =\ell$.  That is, in all three cases, $\HH_Z(i,j)$
equals the given upper bound.  In particular, 
$\HH_Z(i,j) = \deg Z$ if $i+j = \ell$ and $i \geq m_1-1$ and
$j \geq m_1-1$.
\end{proof}

\begin{remark}
Note that we have in fact proved a stronger result.  In conjunction
with Lemma \ref{eventualgrowth}, we can describe $\HH_Z(i,j)$ for
all $(i,j)$ with $i+j \geq \max\{m-1,2m_1-2\}$ directly from the multiplicities of
the points.
\end{remark}

\begin{appendix}
\section{Resolutions and B-regularity}
In this section we prove some modified versions of results in \S 7 of \cite{MS} for finitely generated $\N^k$-graded $R$-modules with $R=[x_{1,0},\ldots,x_{1,n_1},
\ldots,x_{k,0},\ldots,x_{k,n_k}]$ where $\deg x_{i,j} = e_i$,
the $i$th standard basis vector of $\Z^k.$ Here,
we use the notation of \cite{MS}.
\begin{lemma}[Lemma 7.1 in \cite{MS}]\label{appendix lemma}
Let $0 \to M' \to M \to M'' \to 0$ be a short exact sequence of finitely generated $\N^k$-graded $R$-modules.
 If $i \ge 1,$ then 
\[\left( \bigcup _{1\leq j \leq k} (-e_j +\reg^{i+1}(M'))\right)\cap \reg^i(M) \subseteq \reg^i(M'').\]
Otherwise,  when $i = 0$\[\bigcap_{1 \leq j \leq k} (-e_j +\reg^1(M')) \cap \reg^0(M) \subseteq \reg^0(M'').\]
 \end{lemma}

\begin{proof}
For the first statement, the proof proceeds as in \cite{MS}.  The key inclusion
is that for any $j$, $m+\N^k[1-t] \subseteq m+e_j +\N^k[-t]$ for all $t \ge 1.$  However,
the inclusion fails if $t =0$.

Turning to the second part of the claim, we let
\[ m \in \bigcap_{1 \leq j \leq k} (-e_j +\reg^1(M')) \cap \reg^0(M).\]
 The long exact sequence in cohomology gives
\[\cdots \to H^0_B(M)_p \to H^0_B(M'')_p \to H^1_B(M')_p \to \cdots\]
We know that the $H^0_B(M)_p = 0$ for all 
\[p \in m + \N^k[1] = \bigcup_{j = 1}^{k} m+e_j+\N^k.\] 
We will be done if we can show that $H^1_B(M')_p = 0$ for all $p \in m + \N^k[1]$ since 
the vanishings of the higher cohomology modules follow from the first (and stronger) part of the lemma.

We know that $H^1_B(M')_p = 0$ for all $p \in \reg^1(M').$  
So, it is enough to show that $m +\N^k[1] \subseteq \reg^1(M').$   
Since $m \in \bigcap_{1 \leq j \leq k} (-e_j +\reg^1(M')),$ $m +e_j \in \reg^1(M')$ for every $j.$  That is, $m + e_j + \N^k \subseteq \reg^1(M)$ for every $j,$ and since \[\bigcup_{j=1}^{k} m+e_j+\N^k = m+\N^k[1],\] and we are done.
\end{proof}

In our situation, we have the following variant of Corollary 7.3 in \cite{MS}.
\begin{theorem} \label{appendix theorem}
Let $0 \to E_r \to \cdots \rightarrow E_3 \rightarrow E_2 \rightarrow E_1 \rightarrow
E_0 \rightarrow 0$ be a free $\N^k$-graded resolution of a finitely generated 
$\N^k$-graded $R$-module $M.$   Let $m = \min\{r, N+1\}.$  We have
\[\bigcup_{\phi:[m] \rightarrow [k]} \left( \bigcap_{1 \leq i \leq m}
(-e_{\phi(2)} - \cdots - e_{\phi(i)} + \reg^i(E_i))
\right) \subseteq \reg^1(K_0)\]
where $K_0$ is the first syzygy module of $M$ and the union is over all functions $\phi:[m] \rightarrow [k].$ 
\end{theorem}

The result follows from the proof of Theorem 7.2 in \cite{MS} which
is based upon the first conclusion
of Lemma \ref{appendix lemma} and descending induction
on $i$.  
\end{appendix}


\begin{thebibliography}{99}

\bibitem{ACD} A. Aramova, K. Crona, E. De Negri, Bigeneric
initial ideals, diagonal subalgebras and bigraded Hilbert functions.
J. Pure Appl. Algebra {\bf 150} (2000) 215--235.

\bibitem{AH} A. Aramova, J. Herzog, Almost regular sequences and 
Betti numbers.  American J. Math. {\bf 122} (2000) 689--719.

\bibitem{CGG}
M.V. Catalisano, A. V. Geramita, A. Gimigliano, 
Ranks of tensors, secant varieties of Segre varieties and fat points.
Linear Algebra Appl. 
{\bf 355} (2002) 263--285. 

\bibitem{CGG2}
M.V. Catalisano, A. V. Geramita, A. Gimilgliano,
Higer secant varieties of Segre-Veronese varieties.
Projective Varieties with Unexpected Properties: A Volume in
Memory of Giuseppe Veronese (C. Ciliberto, A.V. Geramita, B. Harbourne, R. 
Miro-Roig, K. Ranestad, eds.), de Gruyter, to appear.


\bibitem{CTV} M.V. Catalisano, N.V. Trung, G. Valla,  
A sharp bound for the regularity index of fat points in general position. 
Proc. Amer. Math. Soc. {\bf 118} (1993) 717--724.

\bibitem{Co}
CoCoATeam, CoCoA: a system for doing Computations
 in Commutative Algebra, Available at {\tt http://cocoa.dima.unige.it}

\bibitem{C}
D. Cox, The homogeneous coordinate ring of a toric variety. 
J. Alg. Geom. {\bf 4} (1995) 17--50.

\bibitem{CoxII}
D. Cox, Equations of parametric curves and surfaces via syzygies.  Symbolic Computation:
Solving Equations in Algebra, Geometry and Engineering.  Contemporary Mathematics, vol 286,
AMS, Providence, RI (2001) 1-20.

\bibitem{DG} E.D. Davis, A.V. Geramita,  
The Hilbert function of a special class of $1$-dimensional Cohen-Macaulay graded algebras. 
The curves seminar at Queen's, Vol. III (Kingston, Ont., 1983) Exp. No. H, 29 pp., 
Queen's Papers in Pure and Appl. Math., 67, 
Queen's Univ., Kingston, ON, 1984. 

\bibitem{DS} H. Derksen, J. Sidman, Castelnuovo-Mumford regularity by
  approximation, Adv. Math. {\bf 188} (2004) 104--123.

\bibitem{FL} G. Fatabbi, A. Lorenzini, 
On a sharp bound for the regularity index of any set of fat points. 
J. Pure Appl. Algebra {\bf 161} (2001) 91--111.

\bibitem{EG} D. Eisenbud, S. Goto,  
Linear free resolutions and minimal multiplicity. 
J. Algebra {\bf 88} (1984) 89--133.

\bibitem{GS}
D. Grayson, M. Stillman, Macaulay 2 -- a system for computation in 
algebraic geometry and commutative algebra, http://www.math.uiuc.edu/Macaulay2.


\bibitem{GuVT} E. Guardo, A. Van Tuyl, Fat Points in $\popo$ and their 
Hilbert functions, Can. J. Math. {\bf 56} (2004) 716--741.

\bibitem{ha} H. T. H\`a, Multigraded Castelnuovo-Mumford regularity, $a^*$-invariants and the minimal free resolution.  (2005) Preprint. {\tt math.AC/0501479}.

\bibitem{HVT} H. T. H\`a, A. Van Tuyl, The regularity of 
points in multi-projective spaces.
J. Pure Appl. Algebra {\bf 187} (2004) 153-167. 

\bibitem{H} B. Harbourne, 
Problems and progress: a survey on fat points in $\mathbf P\sp 2$. 
Zero-dimensional schemes and applications (Naples, 2000) 85--132, Queen's
Papers in Pure and Appl. Math., 123, Queen's Univ., Kingston, ON, 2002.

\bibitem{HW} J. W. Hoffman, H. Wang, Castelnuovo-Mumford Regularity in 
Biprojective Spaces, Adv. in Geom. {\bf 4} (2004) 513--536.

\bibitem{MS} D. Maclagan, G. G. Smith, Multigraded Castelnuovo-Mumford 
Regularity, J. Reine Angew. Math. {\bf 571} (2004) 179--212.

\bibitem{R}  T. R\"omer, Homological properties of bigraded regularity.
Ill. J. Math. {\bf 45} (2001) 1361-1376.

\bibitem{SVW} J. Sidman, A. Van Tuyl, H. Wang, Multigraded regularity: coarsenings and resolutions. (2005) Preprint. {\tt math.AC/0505421}.

\bibitem{T} N.V. Trung, 
The Castelnuovo regularity of the Rees algebra and the associated graded ring. 
Trans. Amer. Math. Soc. {\bf 350} (1998) 2813--2832.

\bibitem{TV} N.V. Trung, G. Valla, Upper bounds for the regularity
index of fat points, J. Algebra {\bf 176} (1995) 182-209.

\bibitem{VT} A. Van Tuyl, 
The border of the Hilbert function of a set of points in 
$\pnk$.
J. Pure Appl. Algebra {\bf 176} (2002)  223--247.

\bibitem{VT1} A. Van Tuyl, The Hilbert functions of ACM sets of points
in $\pnk$. J. Algebra {\bf 264} (2003) 420-441.

\bibitem{ZSCC} J. Zheng, T. Sederberg, E. Chionh, D. Cox, Implicitizing 
rational surfaces with base points using the method of moving surfaces. 
Algebraic Geometry and Geometric Modeling, edited by R. Goldman and R. Krausaukas, Contemporary Mathematics, {\bf 334} (2004) 151-168.
\end{thebibliography}
\end{document}